\documentclass[11pt,leqeq]{article}
\usepackage{amssymb,amsthm}
\usepackage{amsmath}
\usepackage{amscd}
\usepackage{xy}
\xyoption{all}
\usepackage[dvips]{graphicx}
\newtheorem{thm}{Theorem}
\newtheorem{cor}[thm]{Corollary}
\newtheorem{prop}[thm]{Proposition} 

\newtheorem{defn}[thm]{Definition}

\newcommand{\des}{\displaystyle}
\date{}

\begin{document}
\setlength{\baselineskip}{16pt}
\title{Compositional Bernoulli numbers}
\author{H\'ector Bland\'\i n and Rafael D\'\i az }
\maketitle

\begin{abstract}
We define and study the combinatorial properties of compositional
Bernoulli numbers and polynomials within the framework of rational
combinatorics.\\

\end{abstract}

\section{Introduction}

Hitherto combinatorial analysis has focus on the study of natural
numbers, yet the time has come to address other numerical
structures by combinatorial means. In
\cite{Blan2} we proposed a framework for the study of the combinatorics of rational
numbers which we review below; first we illustrate how it works
with a simple example. Looking at sequences $1,2,4,8,
\cdots, 2^n, \cdots,$ and $1,1,2,6,24, \cdots, n!,
\cdots,$ combinatorialists will agree that they count the number of subsets and the number of permutations
of a set with $n$ elements, respectively. It is however less clear
what a sequence such as
$$1,2,2,\frac{4}{3},\frac{2}{3},\frac{4}{15},
\cdots,\frac{2^n}{n!}, \cdots$$ might count. The main difference is that whereas
the former are sequences of natural numbers, the later is a
sequence of rational numbers. One can assign combinatorial meaning
to sequences of rational numbers using the notion of cardinality
of groupoids introduced by Baez and Dolan in
\cite{BaezDolan}. In order to find out what the sequence above
"counts" one should find a sequence of groupoids $x_0, x_1,
\cdots, x_n, \cdots$ such that $|x_n|=\frac{2^n}{n!}$. For $n
\geq 0,$ let $x_n$ be the groupoid  whose set of objects $Ob(x_n)$
is the collection of subsets of  $[n]=\{1,2,...,n \}.$ Morphisms
in $x_n$ from $a$ to $b$ are bijections $\alpha:[n]
\longrightarrow [n]$ such that $\alpha(a)=b.$ By definition the cardinality of $x_n$ is
given by
$$|x_n|=\sum_{a \in D(x_n)}\frac{1}{|x_n(a,a)|},$$
where $D(x_n)$ is the set of isomorphisms classes of objects in
$x_n,$ and $x_n(a,a)$ is the set of morphisms from $a$ to $a$ in
$x_n$. In the present case we have that
$$|x_n|=\sum_{k=0}^n \frac{1}{k! (n-k)!}=\frac{1}{n!}\sum_{k=0}^n \left(\begin{array}{c} n\\k
\end{array}\right)=\frac{2^n}{n!},$$ and thus the sequence of
groupoids $x_n$ provides a combinatorial interpretation for the
sequence of rational numbers $\frac{2^n}{n!}$. The goal of
rational combinatorics, from this point of view, is to uncover the
relationship between sequences of rational numbers
and sequences of finite groupoids.\\

We shall actually adopt the functorial viewpoint of Joyal
\cite{j2}, and work with rational species instead of sequences of
groupoids. For a comprehensive study of combinatorial species the
reader may consult \cite{Bergeron}. A presentation of the theory
of species in a categorical context
\cite{SMacLane} with applications to non-commutative spaces is given in
\cite{RDEP}. Rational species were introduced in
\cite{Blan2} and we shall adopt the notation and conventions of
that paper. Further applications of the theory of rational species
are developed in \cite{cas1} and in the forthcoming works
\cite{cas2, DP}. The main ingredients in the definition of the category
of rational species are $\mathbb{B}$ the category of finite sets
and bijections, and $gpd$ the category of finite groupoids. An
object of $gpd$ is a category $G$ such that all its morphisms are
invertible, $Ob(G)$ is a finite set, and  for $x,y \in Ob(G)$ the
set $G(x,y)$ of morphisms from $x$ to $y$ is also finite. Disjoint
union $\sqcup$ and Cartesian product $\times$ give $gdp$ a couple
of monoidal structures with units $\emptyset$, empty groupoid, and
$1$, groupoid with one object and one morphism, respectively.
Cardinality for groupoids yields a valuation map $|\
\ |: Ob(gpd) \longrightarrow
\mathbb{Q}_+$  with values in the semi-ring of
non-negative rational numbers which satisfies: $|x|=|y|$ if $x$
and $y$ are isomorphic, $|x
\sqcup y|= |x|+|y|,$ $|x \times y|=|x||y|,$ $|\emptyset|=0$, and
$|1|=1.$\\

The category of non-negative rational species $gpd^\mathbb{B}$ is
the category of functors from $\mathbb{B}$ to $gpd$; morphisms in
$gpd^\mathbb{B}$ are natural transformations. Monoidal structures
sum and product on $gpd^\mathbb{B}$ are given by
$(F+G)(x)=F(x)\sqcup G(x)$ and
$$(FG)(x)=
\bigsqcup_{x_1\sqcup x_2 =x}F(x_1)
\times G(x_2),$$ where $F$ and $G$ are rational species, and $x, x_1, x_2$
are finite sets. Units for sum and product are $0$ the species
sending a finite set into $\emptyset$, and $1$ the species sending
a non-empty set into $\emptyset$ and the empty set into $1$. The
valuation map $|\ \ |: Ob(gpd^\mathbb{B})
\longrightarrow \mathbb{Q}_+[[x]]$ given by $|F|=
\sum_{n=0}^{\infty}F([n])\frac{x^n}{n!},$ where  $\mathbb{Q}_+[[x]]$ is the semi-ring of formal power
series with coefficients in $\mathbb{Q}_+,$  is such that
$|F|=|G|$ if $F$ and $G$ are isomorphic, $|F+G|=|F| + |G|,$
$|FG|=|F||G|$, $|1|=1$, and $|0|=0.$ The valuation $|F|$ of a
species $F$ is called its generating series. Thus the main problem
of rational combinatorics is: given a non-negative rational
species $F$ find its generating series $|F|
\in \mathbb{Q}_+[[x]]$. We also consider the inverse problem:
given $f \in
\mathbb{Q}_+[[x]]$, find a nice rational species $F$ such that
$|F|=f.$ For example, consider the hyper-exponential power series
$e_{k} \in
\mathbb{Q}_+[[x]]$ given by
$$e_{k}=\sum_{n=0}^{\infty}\frac{x^n}{(n!)^{k}}.$$ Clearly
$e_1$ is the exponential series and, for $k \geq 2$, $e_k$ is a
divided power series with rational coefficients. Let $E_k:
\mathbb{B} \longrightarrow gpd$  be such that for  a finite set$x$, objects in $E_k(x)$ are tuples $(a_1,\cdots, a_{k-1}) \in
P(x)^{k-1}$ where $P(x)$ is the set of subsets of $x$. Morphisms
in $E_k(x)$ from $(a_1,\cdots, a_{k-1})$ to $(b_1,\cdots,
b_{k-1})$ are tuples $(\alpha_1,\cdots,
\alpha_{k-1})$ where $\alpha_i$ is a permutation of $x$ such that
$\alpha_i(a_i)=b_i.$ The generating series of $E_k$ is given by
$$|E_k|=\sum_{n=0}^{\infty}|E_k([n])|\frac{x^n}{n!}=
\sum_{n=0}^{\infty}\left( \sum_{s_1,...,s_{k-1}=1}^n\frac{1}{\prod_i s_i!(n-s_i)!} \right)\frac{x^n}{n!}=
\sum_{n=0}^{\infty}\frac{2^{(k-1)n}}{(n!)^{k-1}}\frac{x^n}{n!}= e_k(2^{k-1}x).$$
Thus $E_k$ provides a combinatorial interpretation for the formal
power series $e_k(2^{k-1}x).$\\

The category of non-negative rational species let us define
combinatorial interpretations for sequences of non-negative
rational numbers. However, both Bernoulli and compositional
Bernoulli numbers require that we are able to consider the
combinatorics of sequences of arbitrary rational numbers,
including negative ones. Let $\mathbb{Z}_2$-$gpd$ be the category
of $\mathbb{Z}_2$-graded finite groupoids, i.e. finite groupoids
$G$ together with a map $Ob(G) \longrightarrow
\mathbb{Z}_2$ sending $x \in Ob(G)$ to $\overline{x} \in
\mathbb{Z}_2$, such that if there is a morphisms in $G$ from $x$
to $y$ then $\overline{x}=\overline{y}.$ Morphisms in
$\mathbb{Z}_2$-$gpd$ are grading preserving morphisms in $gpd$.
$\mathbb{Z}_2$-$gpd$ has monoidal structures disjoint union and
Cartesian product, where the grading on the disjoint union of
groupoids is the disjoint union of the respective gradings, and
the grading on the Cartesian product of groupoids $G$ and $H$ is
$\overline{(x,y)} = \overline{x} + \overline{y}.$  The valuation
map $|\ \ |:\mathbb{Z}_2 \mbox{-}gpd
\longrightarrow \mathbb{Q}$ given by $$|G|= \sum_{x \in D(G)}\frac{(-1)^{\overline{x}}}{|G(x,x)|}$$ is such that
for all $\mathbb{Z}_2$-graded groupoids $G$ and $H$ we have that:
$|G|=|H|$ if $G$ and $H$ are isomorphic, $|G
\sqcup H|= |G|+ |H|$, $|G \times H|= |G| |H|,$ $|\emptyset|=0$, and $|1|=1.$
There is a negative functor $-:\mathbb{Z}_2$-$gpd^{\mathbb{B}}
\longrightarrow
\mathbb{Z}_2$-$gpd^{\mathbb{B}}$ which is the identity both on objects and morphisms but
$\overline{-x}=\overline{x} + 1.$ The category of rational species
$\mathbb{Z}_2$-$gpd^{\mathbb{B}}$ is the category of functors from
$\mathbb{B}$ to $\mathbb{Z}_2$-$gpd$. One defines monoidal
structures sum and product on $\mathbb{Z}_2$-$gpd^{\mathbb{B}}$
and the valuation map $|\ \ |:\mathbb{Z}_2$-$gpd^{\mathbb{B}}
\longrightarrow
\mathbb{Q}[[x]]$ in complete analogy with the case of non-negative rational
species; the resulting structures enjoy similar properties to
those stated for non-negative rational species.\\

The rest of this work is organized as follows. In Section
\ref{cghf} we provide a combinatorial interpretation for Gauss
hypergeometric functions with rational parameters. In Sections
\ref{cbn} and \ref{bnrs} we introduce a generalization of Bernoulli
numbers and provide combinatorial interpretation for such numbers;
we specialize our construction to a variety of interesting
examples. In Section
\ref{cbp} we study the combinatorics of Bernoulli polynomials.
Section \ref{ccbn} contains the main results of this work, namely,
we introduce compositional Bernoulli numbers and provide
combinatorial interpretation for such numbers. In Section
\ref{ccbp} we discuss compositional Bernoulli polynomials.

\section{On the combinatorics of Gauss hypergeometric functions}\label{cghf}

In this section we study the combinatorics of Gauss hypergeometric
functions with rational parameters generalizing the results of
\cite{Blan} where the case of positive rational parameters was tackle.
Consider formal power series of the form
$$h(\pm\frac{a}{b},\pm\frac{c}{d};\pm\frac{e}{f}) = \sum_{n=0}^{\infty}
\frac{(\pm\frac{a}{b})_{n}(\pm\frac{c}{d})_n}{(\pm\frac{e}{f})_{n}} \frac{x^{n}}{n!},$$
where $a,b,c,d,e,f$ are positive natural numbers,  $\pm$ indicates
that a choice of sign has been made, and $(x)_n$ is the Pochhammer
symbol $(x)_n = x(x + 1)(x + 2) . . .(x + n - 1)$ also known as
the increasing factorial \cite{GCRota}. We also need the
Pochhammer $k$-symbol $$(x)_{n,k} = x(x + k)(x + 2k)...(x + (n -
1)k)$$ introduced in
\cite{DP} and further applied in \cite{DP1, DP2, CT}. We proceed to construct functors
$$H(\pm
\frac{a}{b},\pm \frac{c}{d}; \pm \frac{e}{f}):\mathbb{B}
\longrightarrow \mathbb{Z}_2 \mbox{-} gpd,$$ such that
$|H(\pm\frac{a}{b},\pm\frac{c}{d};\pm\frac{e}{f})|=h(\pm\frac{a}{b},\pm\frac{c}{d};\pm\frac{e}{f}).$
The functor $H(\pm\frac{a}{b},\pm\frac{c}{d};\pm\frac{e}{f})$ is
given on a finite set $x$ by
$$H(\pm\frac{a}{b},\pm\frac{c}{d};\pm\frac{e}{f})(x)=
(\pm[a])_{|x|,[b]}(\pm[c])_{|x|,[d]}\overline{\mathbb{Z}}_{b}^{x}
\overline{\mathbb{Z}}_{d}^{x}\overline{\mathbb{Z}}_{\pm e,|x|,f}[f]^{x},$$
where we are making use the following conventions:
\begin{enumerate}
\item A set $x$ may be regarded as the groupoid whose set of objects
is $x$ and whose morphisms are identities.

\item If $G$ is a group then  $\overline{G}$ is the groupoid
with a unique object $1$ and $\overline{G}(1,1)=G.$
\item If $G$ is a groupoid and $x$ a finite set then $G^x$ is the
groupoid whose objects are maps from $x$ into $Ob(G)$, morphisms
in $G^x$ from $f$ to $g$ are given by $G^x(f,g)=\prod_{a \in
Ob(G)}G(f(a),g(a)).$ It is easy to see that $|G^x|=|G|^{|x|}.$

\item If $G$ and $K$ are  groupoids then $( G )_{n,K}=\prod_{i=0}^{n-1}(G
\sqcup (K \times [i])).$ One can show that $|( G )_{n,K}|=(|G|)_{n,|K|}.$

\item  For $n\geq 1,$
$\mathbb{Z}_{n}$ denotes the cyclic group with $n$ elements; we
also set $\overline{\mathbb{Z}}_{-n}=-\overline{\mathbb{Z}}_{n}$.
For $m
\in \mathbb{Z} \backslash \{0 \}$ and $n,l \in  \mathbb{N}$ we set
$\overline{\mathbb{Z}}_{m,n,l}=\prod_{i=0}^{n-1}\overline{\mathbb{Z}}_{m
+ il}.$ One can show that
$|\overline{\mathbb{Z}}_{m,n,l}|=\frac{1}{(m)_{n,l}}.$

\end{enumerate}

\begin{thm}
$$|H(\pm\frac{a}{b},\pm\frac{c}{d};\pm\frac{e}{f})|=h(\pm\frac{a}{b},\pm\frac{c}{d};\pm\frac{e}{f}).$$
\end{thm}
\begin{proof}
$$|(\pm[a])_{n,[b]}(\pm[c])_{n,[d]}\overline{\mathbb{Z}}_{b}^{[n]}
\overline{\mathbb{Z}}_{d}^{[n]}\overline{\mathbb{Z}}_{\pm e,n,f}[f]^{[n]}|= \frac{(\pm a)_{n,b}}{b^n} \frac{(\pm c)_{n,d}}{d^n}
\frac{f^n}{(\pm e)_{n,f}}= \frac{(\pm \frac{a}{b})_{n}(\pm \frac{c}{d})_n}{(\pm\frac{e}{f})_{n}}.$$
\end{proof}
\section{Combinatorics of Bernoulli numbers}\label{cbn}

We introduce a generalization of Bernoulli numbers which may be
motivated as follows. Suppose we are interested in finding a right
inverse for a linear operator $O:V \longrightarrow V$, i.e. an
operator $G:V
\longrightarrow V$ such that $O(G(v))=v$ for $v
\in V$. Assume $O$ can be written as
$$O=f(D)-\pi_{N}(f)(D),$$ where $D: V \longrightarrow V$ is a linear
map for which a right inverse $I: V \longrightarrow V$ is known,
$f \in
\mathbb{Q}[[x]]$ is a formal power series such that $f_N=1$, and for $N \geq 1$, and the $N$-projection map
$$\pi_{N}:\mathbb{Q}[[x]]\rightarrow\mathbb{Q}[[x]]/\left(x^{N}\right)$$
is given by
$$\pi_{N}\left(\sum_{n=0}^{\infty}f_{n}\frac{x^{n}}{n!}\right)=
\sum_{n=0}^{N-1}f_{n}\frac{x^{n}}{n!}.$$
\begin{defn}{\em For
$f=\sum_{n=0}^{\infty}f_{n}\frac{x^n}{n!}\in\mathbb{Q}[[x]]$ and
$N \geq 0$ such that $f_{N}=1$,  the Bernoulli numbers
$B_{N,n}^{f}$ associated with $f$ and $N$ are given by
$$\sum_{n=0}^{\infty}B_{N,n}^{f}\frac{x^{n}}{n!}=
\frac{x^{N}/N!}{f-\pi_{N}f}= \left( N!\sum_{n=0}^{\infty}\frac{n!f_{n+N}}{(n+N)!}\frac{x^n}{n!} \right)^{-1}.$$}
\end{defn}

\begin{thm}\label{ri}
{\em  Under the above conditions a right inverse $G$ for $O$ is
given by
$$G=N! \sum_{n=0}^{\infty}B_{N,n}^{f}\frac{D^{n} \circ I^N}{n!}.$$}
\end{thm}
Classical Bernoulli numbers $B_n=B_{1,n}^{e^x}$ arise when
computing a right inverse $S$ for the finite difference operator
$\Delta(f)=f(x+1) - f(x),$ which acts on functions depending on a
real variable $x$. Let $D=\frac{\partial}{\partial x}$ and $R$ be
a right inverse for $D$, i.e. $R$ is the Riemann integral. Since
$\Delta = e^{D}-1,$ then $S$ is given by
$$S= R + \sum_{n=1}^{\infty}B_{n}\frac{D^{n-1}}{n!}  \ \ \mbox{ where } \ \  \sum_{n=0}^{\infty}B_{n}\frac{x^{n}}{n!}=
\frac{x}{e^x-1}.$$

By definition Bernoulli numbers $B_{N,n}^{f}$ satisfy the
identities
\begin{equation*}
\sum_{k=0}^{n-N}\left(\begin{array}{c} n\\k
\end{array}\right)f_{n-k}B_{N,k}^{f}=\delta_{n,N}\ \ \textrm{for}\ \
n\geq N,
\end{equation*}
and thus can be computed using the recursions:

\begin{equation*}
B_{N,n}^{f}=\left(\begin{array}{c}
N+n \\
n
\end{array} \right)^{-1} \sum_{k=0}^{n-1}\left(
\begin{array}{c} N+n \\ k
\end{array}\right)f_{N+n-k}B_{N,k}^{f}.
\end{equation*}

We provide a combinatorial interpretation for Bernoulli numbers
$B_{N,n}^{f}$ assuming that a combinatorial interpretation for $f$
is known, that is, given $F:\mathbb{B}\longrightarrow
\mathbb{Z}_2$-$gpd$ we construct  $B_N^F:
\mathbb{B}
\longrightarrow \mathbb{Z}_2$-$gpd$ such that $$|B_N^F|=
\frac{x^{N}/N!}{|F|-\pi_{N}|F|}.$$
Consider first the problem of finding a combinatorial
interpretation for $f^{-1}$ assuming that a combinatorial
interpretation for $f$ is known. For $F
\in {\mathbb{Z}_{2}\mbox{-}gpd}^{\mathbb{B}}$ such that
$F(\emptyset)=0,$ let $(1+F)^{-1}:\mathbb{B}\longrightarrow
\mathbb{Z}_{2}\mbox{-}gpd$ send $x \in \mathbb{B}$ into
$$(1+F)^{-1}(x)= \bigsqcup_{x_{1}\sqcup\dots\sqcup x_{k}=x, k \geq 1}(-1)^{k}
\prod_{i=1}^kF(x_i).$$

\begin{prop}\label{Teor29}{\em $(1+F)^{-1}:\mathbb{B}^{n}\longrightarrow
\mathbb{Z}_{2}\mbox{-}gpd$ satisfies $|(1+F)^{-1}|=(1+|F|)^{-1}.$}
\end{prop}
The decreasing factorial rational species
$\mathbb{Z}_{(N)}:\mathbb{B}\longrightarrow gpd$ is such that for
$x \in \mathbb{B}$ we have $Ob(\mathbb{Z}_{(N)}(x))=\{x\}$ if $|x|
\geq N$ and empty otherwise. For $|x|
\geq N$ morphisms in $\mathbb{Z}_{(N)}(x)$ are given by
$$\mathbb{Z}_{(N)}(x,x)=\mathbb{Z}_{|x|}\times\mathbb{Z}_{|x|-1}\times\dots \times \mathbb{Z}_{|x|-N+1}.$$
The derivative $\partial F$ of a species $F$ is given by $\partial
F(x)=F(x \sqcup \{x\})$ for $x \in \mathbb{B}$; it is easy to
check that $|\partial F|= \partial |F|$. For $n \geq 1$ and $G$ a
groupoid,  let $nG$ be the groupoid $G\sqcup ...
\sqcup G$ ($n$ copies).
\begin{prop}{\em
Let $F:\mathbb{B}\longrightarrow \mathbb{Z}_{2}\mbox{-}gpd$ be a
rational species such that $F([N])=1$, then
\begin{equation*}
\left|1 + N!\partial^{N}\left(F\times
\mathbb{Z}_{(N)}\right)\right|=N!\frac{|F|-\pi_{N}|F|}{x^{N}}.
\end{equation*}}
\end{prop}
\begin{proof} For $n=0$ the desired result is obvious. For $n \geq 1$ we have
$$|1 + N!\partial^{N}\left(F\times
\mathbb{Z}_{(N)}\right)([n])|=N!|F([n] \sqcup [N])||\mathbb{Z}_{(N)}([n] \sqcup [N])|=N!|F|_{n+N}\frac{n!}{(n+N)!},$$
which is the $n$-th coefficient of the divided power series
$N!\frac{|F|-\pi_{N}(|F|)}{x^{N}}.$
\end{proof}

\begin{thm}\label{aaa}{\em
Let $F:\mathbb{B}\longrightarrow \mathbb{Z}_{2}\mbox{-}gpd$ be
such that $F\left([n]\right)=1$. The species $B_N^F:\mathbb{B}
\longrightarrow \mathbb{Z}_{2}\mbox{-}gpd$ sending $x \in
\mathbb{B}$ into $$B_N^F(x) =\bigsqcup_{x_{1}\sqcup\dots\sqcup x_{k}=x, k
\geq 1}(-N!)^k
\prod_{i=1}^{k}F \times \mathbb{Z}_{(N)}(x_i \sqcup [N]),$$
is such that $\left|B_N^F\right| =
\sum_{n=0}^{\infty}B_{N,n}^{|F|}\frac{x^{n}}{n!}.$}
\end{thm}
\begin{proof} Follows from Proposition \ref{Teor29} and the identity $|B_N^F| =
|1 + N!\partial^{N}\left(F\times
\mathbb{Z}_{(N)}\right)|^{-1}.$
\end{proof}

Next we compare  $\Delta$ with $D=\frac{\partial}{\partial x}$;
for example we may like to know a right inverse for the operator
$\Delta - D=e^D -1 -D$, or more generally a right inverse for the
operator $e^D -
\pi_N(e^D).$ According to Theorem
\ref{ri} a right inverse for $e^D - \pi_N(e^D)$ is given by
$$G=N! \sum_{n=1}^{N-1}B_{N,n}\frac{I^{N-n}}{n!} +
N! \sum_{n=N}^{\infty}B_{N,n}\frac{D^{n-N}}{n!},$$ where the
Bernoulli  numbers $B_{N,n}$ are such that
\begin{equation*}
\frac{x^{N}/N!}{e^{x}-1-x}=\sum_{n=0}^{\infty}B_{N,n}\frac{x^{n}}{n!}.
\end{equation*}
Explicitly the first Bernoulli numbers $B_{N,n}=B_n$ are shown in
the table:
\begin{center}
\footnotesize{
\begin{tabular}{|c|c|c|c|c|c|c|c|c|c|c|c|c|c|c|c|}
\hline
      $n$ & 0 & 1 & 2 & 3 & 4 & 5 & 6 & 7 & 8 & 9 & 10 & 11 & 12 & 13 & 14 \\
\hline
  $B_{n}$ & 1 & -1/2 & 1/6 & 0 & -1/30 & 0 & 1/42 & 0 & -1/30 & 0 & 5/66 & 0 & -691/2730 & 0 &7/6\\
\hline
\end{tabular}
}
\end{center}
The first Bernoulli numbers $B_{2,n}$ are given by:
\begin{center}
\footnotesize{
\begin{tabular}{|c|c|c|c|c|c|c|c|c|c|c|}
\hline
  $n$ & 0 & 1 & 2 & 3 & 4 & 5 & 6 & 7 & 8 & 9 \\
\hline
  $B_{2,n}$ & 1 & -1/3 & 1/18 & 1/90 & -1/270 & -5/1134 & -1/5670 & 7/2430 & 13/7290 & -307/133650 \\
\hline
\end{tabular}}
\end{center}
\bigskip

The exponential species $E$ given by $E(x)=\{x \}$ for $x \in
\mathbb{B}$ is such that $|E|= e ^x.$  Theorem \ref{aaa} implies that the
generating series of the species $B_N:\mathbb{B} \longrightarrow
\mathbb{Z}_{2}\mbox{-}gpd$ sending $x \in \mathbb{B}$ into
$$B_N(x)= \bigsqcup_{x_1 \sqcup ... \sqcup x_k=x, k \geq 1}(-N!)^k
\prod_{i=1}^{k}\mathbb{Z}_{(N)}(x_i \sqcup [N])$$
is $\sum_{n=0}^{\infty}B_{N,n}^{f}\frac{x^{n}}{n!}$. Thus we have
obtained a combinatorial interpretation for $B_{N,n}$ in terms of
the cardinality of $\mathbb{Z}_{2}$-graded groupoids:
\begin{eqnarray*}
B_{N,n}= |B_N([n])|=n!
\des{\sum_{a_{1}+\dots+a_{k}=n, k \geq 1}\frac{(-N!)^k}{(a_{1}+N)!\dots(a_{k}+N)!}}.
\end{eqnarray*}

Let us consider Bernoulli numbers associated with the sine
function and $N=2L+1$  an odd number. For $N=1$ we obtain
Bernoulli numbers $B_{1,n}^{sin}$ given by
\begin{equation*}
\frac{x}{sin(x)}=\sum_{n=0}^{\infty}B_{1,n}^{sin}\frac{x^{n}}{n!}.
\end{equation*}
The first few values of the sequence $B_{1,n}^{sin}$ are shown in
the table:

\begin{center}
\footnotesize{
\begin{tabular}{|c|c|c|c|c|c|c|c|c|c|c|c|c|c|c|}
\hline
      $n$ & 0 & 1 & 2 & 3 & 4 & 5 & 6 & 7 & 8 & 9 & 10 & 11 & 12 & 13 \\
\hline
  $B_{1,n}^{\ \sin}$ & 1 & 0 & 1/3 & 0 & 7/15 & 0 & 31/21 & 0 & 127/15 & 0 & 2555/33 & 0 & 1414477/1365 & 0 \\
\hline
\end{tabular}
}
\end{center}
One can show that
$B_{1,2n}^{sin}=(-1)^{n-1}\left(2^{2n}-2\right)B_{2n}$ and
$B_{1,2n+1}^{sin}=0$. For $N=3$ we obtain the Bernoulli numbers
$B_{2,n}^{sin}$ given by
\begin{equation*}
\frac{x^{3}/3!}{sin(x)-x}=\sum_{n=0}^{\infty}B_{3,n}^{sin}\frac{x^{n}}{n!}.
\end{equation*}
Explicitly the first few values of the sequence $B_{3,n}^{sin}$
are given in the table:

\begin{center}
\footnotesize{
\begin{tabular}{|c|c|c|c|c|c|c|c|c|c|c|c|c|}
\hline
      $n$ & 0 & 1 & 2 & 3 & 4 & 5 & 6 & 7 & 8 & 9 & 10 & 11 \\
\hline
  $B_{3,n}^{sin}$ & -1 & 0 & -1/10 & 0 & -11/350 & 0 & -17/1050 & 0 & -563/57750 & 0 & -381/250250 & 0 \\
\hline
\end{tabular}
}
\end{center}
Let $Sin$ be the species such that
$$Sin(x)=\left\{\begin{array}{cc} (-1)^n & \mbox{ if \ }
|x|=2n+1 \\ 0 & \mbox{if \ }  x \mbox{ is even }
\end{array}\right.$$
Clearly $|Sin|=sin,$  $B_{2L+1}^{Sin}(x)=0$ if $|x|$ is odd, and
if $|x|$ is even then
\begin{equation*}
B_{2L+1}^{Sin}(x) = \bigsqcup_{x_1 \sqcup ... \sqcup x_k=x, k
\geq 1}(-1)^{\frac{|x|}{2}+kL}(-(2L+1)!)^{k}
\prod_{i=1}^{k}\mathbb{Z}_{(N)}(x_i \sqcup [2L+1]),
\end{equation*}
where the cardinality of each set $x_i$ is even. Therefore we
obtain that
\begin{eqnarray*}
B_{2L+1,2n}^{sin} = |B_{2L+1}^{Sin}[2n]| = (-1)^n 2n!
\des{\sum_{2a_{1}+\dots+2a_{k}=2n, k \geq 1}\frac{(-1)^{k(L+1)}(2L+1)!^{k}}{(2a_{1}+2L+1)!\dots(2a_{k}+2L+1)!}}.
\end{eqnarray*}
Similarly one can consider Bernoulli numbers $B_{N,n}^{cos}$
associated with the cosine function for $N= 2L$ an even number.
For $N=2$ the Bernoulli numbers $B_{2,n}^{cos}$ are such that:
\begin{equation*}
\frac{x^{2}/2!}{cos(x)-1}=\sum_{n=0}^{\infty}B_{2,n}^{cos}\frac{x^{n}}{n!}.
\end{equation*}
The first few values of $B_{2,n}^{cos}$  are shown in the table:
\begin{center}
\footnotesize{
\begin{tabular}{|c|c|c|c|c|c|c|c|c|c|c|c|c|c|c|}
\hline
      $n$ & 0 & 1 & 2 & 3 & 4 & 5 & 6 & 7 & 8 & 9 & 10 & 11 & 12 & 13 \\
\hline
  $B_{2,n}^{\ cos}$ & -1 & 0 & -1/6 & 0 & -1/10 & 0 & -5/42 & 0 & -7/30 & 0 & -15/22 & 0 & -7601/2730 & 0 \\
\hline
\end{tabular}
}
\end{center}
Let $Cos$ be the species such that
$$Cos(x)=\left\{\begin{array}{cc} (-1)^n & \mbox{ if \ }
|x|=2n \\ 0 &  \mbox{ if } x \mbox{ is odd }
\end{array}\right.$$
Clearly $|Cos|=cos$, $B_{2L}^{Cos}(x)=0$ if $|x|$ is odd, and if
$|x|$ is even then
\begin{equation*}
B_{2L}^{Cos}(x) = \bigsqcup_{x_1 \sqcup...\sqcup x_k=x, k
\geq 1}(-1)^{\frac{|x|}{2}+kL}(-2L!)^k
\prod_{i=1}^{k}\mathbb{Z}_{(N)}(x_i \sqcup [N]),
\end{equation*}
where the cardinality of each set $x_i$ is even. Therefore we
obtain that
\begin{eqnarray*}
B_{2L,2n}^{cos} = |B_{2L}^{Cos}[2n]|= (-1)^{n}2n!
\des{\sum_{2a_{1}+\dots+2a_{k}=2n, k \geq 1}\frac{(-1)^{k(L+1)}2L!^k}{(2a_{1}+2L)!\dots(2a_{k}+2L)!}}.
\end{eqnarray*}

\section{Bernoulli numbers for rational species}\label{bnrs}

In this section we consider Bernoulli numbers associated with
formal power series with rational coefficients. Let $M$ be a
positive integer and $\mathbb{Z}^{M}:\mathbb{B}\longrightarrow
gpd$ be the rational species such that for $x \in \mathbb{B}$ the
groupoid $\mathbb{Z}^{M}(x)$ is given by
$$Ob(\mathbb{Z}^{M}(x))=\left\{\begin{array}{cc} \{x\} & \mbox{ if \ }
x\neq\emptyset \\ \emptyset  & \mbox{ if \ } x=\emptyset
\end{array}\right.$$ and  $\mathbb{Z}^{M}(x)(x,x)=\mathbb{Z}_{|x|}^{M}$
for $x$ nonempty; clearly
$\left|\mathbb{Z}^{M}\right|=\sum_{n=1}^{\infty}\frac{1}{n^{M}}\frac{x^{n}}{n!}.$
The first Bernoulli numbers for $\left|\mathbb{Z}\right|$ are
given in the table:

\begin{center}
\begin{tabular}{|c|c|c|c|c|c|c|c|c|c|}
\hline
      $n$ & 0 & 1 & 2 & 3 & 4 & 5 & 6 & 7 \\
\hline
  $B_{1,n}^{\mathbb{Z}}$ & 1 & -1/4 & 1/72 & 1/96 & 61/21600 & -1/640 & -12491/5080320 & -479/580608 \\
\hline
\end{tabular}
\end{center}
For $M=2$ we get the species
$\mathbb{Z}^{2}:\mathbb{B}\longrightarrow gpd$ with generating
series $\left|\mathbb{Z}^{2}\right|=
\displaystyle{\sum_{n=1}^{\infty}\frac{1}{n^{2}}\frac{x^{n}}{n!}}.$
The first Bernoulli numbers for  $\left|\mathbb{Z}^{2}\right|$ are
shown in the table:
\begin{center}
\begin{tabular}{|c|c|c|c|c|c|c|c|c|}
\hline
      $n$ & 0 & 1 & 2 & 3 & 4 & 5 & 6 \\
\hline
  $B_{1,n}^{\mathbb{Z}^{2}}$ & 1 & -1/8 & 11/432 & 1/144 & -217/324000 & -157/64800 & -21503/16669800 \\
\hline
\end{tabular}
\end{center}
According to Theorem \ref{aaa} the generating series of the
species $B_N^{\mathbb{Z}^M}$ sending $x \in
\mathbb{B}$ by
\begin{equation*}
B_N^{\mathbb{Z}^M}(x) = \bigsqcup_{x_1 \sqcup ... \sqcup x_k=x, k
\geq 1}(-N!)^k
\prod_{i=1}^{k}\mathbb{Z}_{(N)}(x_i \sqcup [N])
\prod_{i=1}^{k}\overline{\mathbb{Z}}_{|x_i| + N}^M
\end{equation*}
is $\sum_{n=0}^{\infty}B_{N,n}^{\mathbb{Z}^M}\frac{x^n}{n!}.$
Therefore we obtain that

$$B_{N,n}^{\mathbb{Z}^M}=|B_N^{\mathbb{Z}^M}([n])|= (-1)^{n}n!
\des{ \sum_{a_{1}+\dots+a_{k}=n, k \geq
1}\frac{N!^k}{(a_{1}+N)!(a_{1}+N)^M\dots(a_{k}+N)!(a_{k}+N)^M}}.$$
Let $\mathbb{Z}^{(M)}:\mathbb{B}\longrightarrow gpd$ be the
rational species sending $x \in \mathbb{B}$ into
$\mathbb{Z}^{(M)}(x)$ the groupoid given by
\begin{equation*}
Ob(\mathbb{Z}^{(M)}(x))=\left\{\begin{array}{cc}\{x\}& \mbox{if \
} x\neq\emptyset\\ \emptyset & \mbox{ if \ } x=\emptyset,
\end{array}\right.
\end{equation*}
and
$\mathbb{Z}^{(M)}(x)(x,x)=\mathbb{Z}_{|x|}\times\mathbb{Z}_{|x|+1}
\times\dots\times\mathbb{Z}_{|x|+M-1}$ for $x$ non-empty. The
generating series of $\mathbb{Z}^{(M)}$ is
\begin{equation*}
\left|\mathbb{Z}^{(M)}\right|=\sum_{n=1}^{\infty}\frac{1}{(n)^{(M)}}\frac{x^{n}}{n!}.
\end{equation*}
The first Bernoulli numbers for $|\mathbb{Z}^{(3)}|$ and $N=3$ are
shown in the table:
\begin{center}
\footnotesize{
\begin{tabular}{|c|c|c|c|c|c|c|c|c|c|c}
\hline
      $n$ & 0 & 1 & 2 & 3 & 4 & 5 & 6 & 7 \\
\hline
  $B_{3,n}^{\mathbb{Z}^{3}}$ & 60 & -15/12 & 9/56 & 3/64 & 401/31360 & 127/50176 & -9089/33116160
  & -192233/264929280 \\
\hline
\end{tabular}
}
\end{center}

\section{Generalized Bernoulli polynomials}\label{cbp}

Let $f\in\mathbb{Q}[[x]]$ and $N \geq 1$ be such that $f_{N}=1.$
Bernoulli polynomials $B_{N,n}^{f}(x)$ are such that
\begin{equation*}
\des{\sum_{n=0}^{\infty}B_{N,n}^{f}(x)\frac{y^{n}}{n!}=
\frac{f(xy)\left(y^{N}/N!\right)}{f(y)-\pi_{N}(f)(y)}}.
\end{equation*}
For example Bernoulli polynomials $B_{n}(x)$ and $B_{2,n}(x)$ are
given by the identities:
$$\des{\sum_{n=0}^{\infty}B_{n}(x)\frac{y^{n}}{n!}=\frac{e^{xy}y}{e^{y}-1}}
\mbox{ \  \  and  \  \  }
\des{\sum_{n=0}^{\infty}B_{2,n}(x)\frac{y^{n}}{n!}=
\frac{e^{xy}\left(y^{2}/2!\right)}{e^{y}-1-y}}.$$
It is not hard to check that Bernoulli polynomials
$B_{N,n}^{f}(x)$ satisfy the recursion:
\begin{equation*}
\sum_{k=0}^{n-N}\left(\begin{array}{c} n\\k
\end{array}\right)f_{n-k}B_{N,k}^{f}(x)=\left(\begin{array}{c} n\\N
\end{array}\right)f_{n-N}x^{n-N}\ \textrm{for}\ n\geq N,
\end{equation*}
or, equivalently,  $B_{n,0}^{f}(x)=f_{0}$ and
\begin{equation*}
B_{N,n}^{f}(x)=\frac{\des{\left(\begin{array}{c} n+N\\
N\end{array}\right)f_{n}x^n-\sum_{k=0}^{n-1}\left(\begin{array}{c}
n+N\\k
\end{array}\right)f_{n+N-k}B_{N,k}^{f}(x)}}{\left(\begin{array}{c}n+N \\ n
\end{array}\right)}.
\end{equation*}
Next result writes Bernoulli polynomials $B_{N,n}^{f}(x)$ in terms
of Bernoulli numbers $B_{N,n}^{f}$.
\begin{thm}{\em
Let $f\in\mathbb{Q}[[x]]$ and $N \geq 1$ be such that $f_{N}=1,$
then $B_{N,n}^{f}(x)$ is given by}
\begin{equation*}
B_{N,n}^{f}(x)=\sum_{k=0}^{n}\left(\begin{array}{c} n \\ k
\end{array}\right)B_{N,n-k}^{f}f_{k}x^{k}.
\end{equation*}
\end{thm}
\begin{proof}
$$\sum_{n=0}^{\infty}B_{N,n}^{f}(x)\frac{z^{n}}{n!} = \left(\sum_{n=0}^{\infty}B_{N,n}^{f}\frac{z^{n}}{n!}\right)
\left(\sum_{n=0}^{\infty}f_{n}x^{n}\frac{z^{n}}{n!}\right) =\sum_{n=0}^{\infty}\left(\sum_{k=0}^{n}\left(\begin{array}{c} n
\\ k
\end{array}\right)B_{N,n-k}^{f}f_{k}x^{k}\right)\frac{z^{n}}{n!}.$$
\end{proof}

Consider the species $\mathrm{S}:\mathbb{B}\longrightarrow gpd$
such that $\mathrm{S}(x)$ is the groupoid given by
$Ob\left(\mathrm{S}(x)\right)=\{x\}$ and
$\mathrm{S}(x)\left(x,x\right)=S_{|x|}$ for $x \in \mathbb{B}$.
The generating series of $S$ is
$|\mathrm{S}|=\displaystyle{\sum_{n=0}^{\infty}\frac{1}{n!}\frac{x^{n}}{n!}},$
and the corresponding Bernoulli polynomials are given by
\begin{eqnarray*}
B_{1,0}^{S}(x)&=&1, \ \ B_{1,1}^{S}(x)=x-\frac{1}{4}, \ \ B_{1,2}^{S}(x)=\frac{1}{2}x^{2}-\frac{1}{2}x+\frac{5}{72},\\
B_{1,3}^{S}(x)&=&\frac{1}{6}x^{3}-\frac{3}{8}x^{2}+\frac{5}{24}x-\frac{1}{48},\\
B_{1,4}^{S}(x)&=&\frac{1}{24}x^{4}-\frac{1}{6}x^{3}+\frac{5}{24}x^{2}-\frac{1}{12}x+\frac{139}{21600}, \\
B_{1,5}^{S}(x)&=&\frac{1}{120}x^{5}-\frac{5}{96}x^{4}+\frac{25}{216}x^{3}-\frac{5}{48}x^{2}+\frac{139}{4320}x-\frac{1}{540}.
\end{eqnarray*}

Bernoulli numbers for $S$ are shown in the table:

\begin{center}
\footnotesize{
\begin{tabular}{|c|c|c|c|c|c|c|c|c|c|}
\hline
      $n$ & 0 & 1 & 2 & 3 & 4 & 5 & 6 & 7 & 8  \\
\hline
  $B_{1,n}^S$ & -1 & -1/4 & 5/72 & 1/48 & 139/21600 & -1/540 & 859/2540160 &
   71/483840 & -9769/36288000 \\
\hline
\end{tabular}
}
\end{center}
\section{Compositional Bernoulli numbers}\label{ccbn}

Let $f^{<-1>}$ be the compositional inverse of the formal power
series $f$, i.e. $f\circ f^{<-1>} = x = f^{<-1>} \circ f.$ Also if
$k$ is a positive integer we define inductively $f^{<1>}=f$ and
$f^{<k+1>}=f^{<k>}
\circ f.$
\begin{defn}\label{listo}{\em
Let $f=\sum_{n=0}^{\infty}f_{n}\frac{x^n}{n!}\in\mathbb{Q}[[x]]$
and $N \geq 1$ be such that $f_{N}=1$. The compositional Bernoulli
numbers $C_{N,n}^{f}$ are given by
\begin{equation*}
\des{\sum_{n=1}^{\infty}C_{N,n}^{f}\frac{x^{n}}{n!} = \left(N!x^{1-N}(f-\pi_{N}f)\right)^{<-1>}=
\left(\sum_{n=1}^{\infty}
{\frac{N!n!f_{N+n-1}}{(N+n-1)!}}\frac{x^{n}}{n!}\right)^{<-1>}}.
\end{equation*}}
\end{defn}

In order to compute compositional Bernoulli numbers we use the
recursion below which is a direct consequence of  Definition
\ref{listo}.

\begin{prop}{\em For $n \geq 1$,
\begin{equation*}
C_{N,n}^{f}= \sum_{a_{1}+a_{2}+\dots+a_{k}=n}\frac{(-N!n!)^k}{k!}
\left(\begin{array}{cc} n\\
a_{1},a_{2},\dots,a_{k}\end{array}\right)C_{N,k}^{f}{\prod_{j=1}^{k}\frac{f_{N+
a_{j}-1}} {(N+a_{j}-1)!},}
\end{equation*}
 where $k \geq 2$.}
\end{prop}

Let us provide a combinatorial interpretation for  compositional
Bernoulli numbers $C_{N,n}^{f}$ assuming that a combinatorial
interpretation for $f$ is known, i.e. given  $F:
\mathbb{B}\longrightarrow \mathbb{Z}_2$-$gpd$ we construct $C_N^F: \mathbb{B}
\longrightarrow \mathbb{Z}_2$-$gpd$ such that $$|C_N^F|=
\left(N!x^{1-N}(|F|-\pi_{N}(|F|)\right)^{<-1>}.$$
First consider the problem of finding a combinatorial
interpretation for $f^{<-1>}$ assuming that a combinatorial
interpretation for $f$ is known. Let $Par:\mathbb{B}
\longrightarrow \mathbb{B}$ be the species sending a finite set
$x$ into $Par(x)$ the set of partitions of $x$, i.e. an element
$\pi$ in $Par(x)$ is a family of non-empty subsets of $x$ such
that $\cup_{b \in \pi}b=x$ and $b \cap c =
\emptyset$ for $b,c \in \pi$. We write $a \vdash n$ if $n$ is a positive integer and
$a=(a_1,...,a_l)$ is a sequence of positive integers  such that
$|a|=a_1 + ... + a_l = n.$  The integer $l(a)$ is called the
length of $a$. The generating function of $Par$ is given by
$$|Par|=\sum_{n=1}^{\infty}
\left(\sum_{a\vdash n}\frac{1}{l(a)!}\binom{n}{a_1,...,a_l} \right)\frac{x^n}{n!} .$$
For $d \geq 1$ consider the species $Par_d^s:\mathbb{B}
\longrightarrow \mathbb{B}$ sending $x \in \mathbb{B}$ into
$Par_d^s(x)$ the set of $d$-tuples $\pi = (\pi_1,...,\pi_d)$ such
that: $\pi_1 \in Par(x),$  $\pi_i \in Par(\pi_{i-1})$ for $2
\leq i \leq d,$ $|\pi_d| \geq 2,$  and $|b| \geq s$ for $b \in \pi_i$ and $1 \leq i \leq d.$

\begin{prop}{\em The generating series of $Par_d^s$ is given by
$$|Par_d^s| = \sum_{n=1}^{\infty}
\left( \sum_{a_1,...,a_d}\frac{1}{l(a_d)!}
\prod_{i=1}^d \binom{l(a_{i-1})}{a_{i1},...,a_{il(a_i)}}\right)\frac{x^n}{n!} ,$$ where $a_1 \vdash
l(a_{0}):=n$, $a_i \vdash l(a_{i-1})$ for $2 \leq i \leq d,$
$|l(a_d)| \geq s,$ and $a_{ij} \geq s$. }
\end{prop}

Let $gpd^\mathbb{B}_1$ be the full subcategory of rational species
such that $F(\emptyset)=\emptyset$ and $F([1])=1$.
\begin{prop}{\em If $F$ belongs to $gpd^\mathbb{B}_1$ and $d \geq 1$, then the rational species $F^{<d+1>}$ is given by
$$F^{<d+1>}(x)= \bigsqcup_{\pi \in Par_d^1(x)} F(\pi_d) \times \prod_{i=1}^d \prod_{b \in \pi_i}F(b).$$}
\end{prop}
\begin{proof}
For $d=1$ the formula above is  the well-known result:
$$F^{<2>}(x)= \bigsqcup_{\pi \in Par(x)}  F(\pi)\times \prod_{b \in \pi}F(b).$$
The desired formula follows by induction: $$F^{<d+1>}(x)= F^{<d>}
\circ F = \bigsqcup_{\pi_1 \in Par(x)}
\left( F^{<d>}(\pi_1)\times
\prod_{b \in \pi_1}F(b) \right).$$

\end{proof}

\begin{cor}{\em  Let $F$ be in $gpd^\mathbb{B}_1$ and $d \geq 1$,
then we have that
$$|F^{<d+1>}| = \sum_{n=0}^{\infty}\left( \sum_{a_1,...,a_d}
\frac{|F|_{l(a_d)}}{l(a_d)!} \prod_{i=1}^d \binom{l(a_{i-1})}{a_{i1},...,a_{il(a_i)}}\prod_{j=1}^{l(a_i)}  |F|_{a_{ij}}
\right) \frac{x^n}{n!},$$
where $a_1 \vdash l(a_{0}):=n$ and $a_i \vdash l(a_{i-1})$ for $2
\leq i
\leq d.$}
\end{cor}

For $F$ in ${\mathbb{Z}_{2}\mbox{-}gpd}^{\mathbb{B}}$ such that
$F(\emptyset)=F([1])=0$ consider the species $(X+F)^{<-1>}$ in
${\mathbb{Z}_{2}\mbox{-}gpd}^{\mathbb{B}}_1$ sending $x \in
\mathbb{B}$ with $|x| \geq 2$ into
$$(X+F)^{<-1>}(x)= -F(x) \sqcup
\bigsqcup_{\pi \in Par_d^2(x), d \geq 1}
(-1)^{d+1}  F(\pi_d) \times \prod_{i=1}^d \prod_{b
\in \pi_i}F(b).$$ The disjoint union above is
finite since it is restricted to partitions whose blocks are at
least of cardinality two.

\begin{thm}\label{Teor50}{\em $(X+F)^{<-1>}$ is such that
$|(X+F)^{<-1>}|=(X+|F|)^{<-1>}.$ }
\end{thm}

\begin{proof}
The result follows from the identities:
$$(X+F)^{<-1>}= X - F +
\sum_{d=1}^{\infty}(-1)^{d+1}F^{<d+1>},$$
$$F^{<d+1>}(x)= \bigsqcup_{\pi  \in Par_d^2(x) }F(\pi_d) \times \prod_{i=1}^d\prod_{b \in \pi_i}F(b).$$
\end{proof}

\begin{cor}{\em Let $f=\sum_{n=0}^{\infty}f_{n}\frac{x^n}{n!}\in\mathbb{Q}[[x]]$ be such that  $f_0 = f_1 =0$, then
$$(x +
f)^{<-1>}= x +
\sum_{n=2}^{\infty}\left(-f_n + \sum_{a_1,...,a_d, d \geq 1}(-1)^{d+1}
\frac{f_{l(a_d)}}{l(a_d)!} \prod_{i=1}^d \binom{l(a_{i-1})}{a_{i1},...,a_{il(a_i)}}\prod_{j=1}^{l(a_i)}f_{a_{ij}}\right)\frac{x^n}{n!},$$
where $a_1 \vdash l(a_{0}):=n,$ $a_i \vdash l(a_{i-1})$ for $2
\leq i \leq d,$ $l(a_i) \geq 2$ and $a_{ij} \geq 2$.}
\end{cor}

\begin{prop}\label{pp}{\em
Let $F:\mathbb{B}\longrightarrow \mathbb{Z}_{2}\mbox{-}gpd$ be a
rational species such that $F([N])=1$, then
\begin{equation*}
\left|X + N!\partial^{N-1}\left(F\times
\mathbb{Z}_{(N-1)}\right)\right|=\left(N!x^{1-N}(|F|-\pi_{N}|F|)\right)^{<-1>}.
\end{equation*}}
\end{prop}

From Theorem \ref{Teor50} and Proposition \ref{pp} we obtain the
promised combinatorial interpretation for the compositional
Bernoulli numbers.

\begin{thm}\label{bernoulligeneralizado}{\em
For $F:\mathbb{B}\longrightarrow \mathbb{Z}_{2}\mbox{-}gpd$  such
that $F\left([N]\right)=1$, let $C_N^F \in
{\mathbb{Z}_{2}\mbox{-}gpd}^{\mathbb{B}}_1$ send $x
\in \mathbb{B}$ with $|x| \geq 2$ into $C_N^F(x)$, the disjoint union in $\mathbb{Z}_{2}\mbox{-}gpd$
of $-N! F\times \mathbb{Z}_{(N-1)}(x \sqcup [N-1])$ and
$$\bigsqcup_{\pi \in Par_d^2(x), d \geq 1}(-N!)^{d+1} F\times
\mathbb{Z}_{(N-1)}(\pi_d \sqcup [N-1]) \times \prod_{i=1}^d \prod_{b
\in \pi_i}F \times
\mathbb{Z}_{(N-1)}(b \sqcup [N-1]).$$
Then $|C_N^F| =\sum_{n=0}^{\infty}C_{N,n}^{|F|}\frac{x^{n}}{n!}.$}
\end{thm}

Recall that we denote by $E:\mathbb{B}\longrightarrow
\mathbb{Z}_{2}\mbox{-}gpd$  the exponential species.

\begin{cor}
{\em Let $C_N^{E}:\mathbb{B} \longrightarrow
{\mathbb{Z}_{2}\mbox{-}gpd}^{\mathbb{B}}_1$ be the species sending
$x \in \mathbb{B}$ with $|x| \geq 2$ into $C_N^{E}(x)$, the
disjoint union of the $\mathbb{Z}_{2}$-graded groupoids $-N!
F\times
\mathbb{Z}_{(N-1)}(x \sqcup [N-1])$ and
$$\bigsqcup_{\pi \in Par_d^2(x), d \geq 1} (-N!)^{d+1}
\mathbb{Z}_{(N-1)}(\pi_d \sqcup [N-1]) \times \prod_{i=1}^d \prod_{b
\in \pi_i}\mathbb{Z}_{(N-1)}(b \sqcup [N-1]).$$
Then $|C_N^{E}|
=\sum_{n=0}^{\infty}C_{N,n}^{|E|}\frac{x^{n}}{n!}.$}
\end{cor}

Compositional Bernoulli numbers for $f(x)=e^{x}$ and $N=1$ are
such that
\begin{equation*}
\sum_{n=1}^{\infty}C_{1,n}
\frac{x^{n}}{n!}=\left(e^{x}-1\right)^{<-1>}=\ln(1+x)=\sum_{n=1}^{\infty}(-1)^{n-1}(n-1)!\frac{x^{n}}{n!},
\end{equation*}
thus $C_{1,n}=(-1)^{n-1}(n-1)!$  for $n\geq 1$.  For $f(x)=e^{x}$
and $N=2$ we obtain compositional Bernoulli numbers $C_{2,n}$
which are such that
\begin{equation*}
\sum_{n=1}^{\infty}C_{2,n}
\frac{x^{n}}{n!}=\left(2!\frac{(e^{x}-1-x)}{x}\right)^{<-1>}
=\left(\sum_{n=1}^{\infty}\frac{2}{n+1}\frac{x^n}{n!}
\right)^{<-1>}.
\end{equation*}
The first compositional Bernoulli numbers $C_{2,n}$ are shown in
the table:
\begin{center}
\footnotesize{
\begin{tabular}{|c|c|c|c|c|c|c|c|c|c|c|c|c|}
\hline
      $n$ & 0 & 1 & 2 & 3 & 4 & 5 & 6 & 7 & 8 & 9  \\
\hline
  $C_{2,n}$ & 0 & 1 & -2/3 & 5/6 & -68/45 & 193/54 & -655/53 & 19349/540 &
  -57736/405 & 520343/810\\
\hline
\end{tabular}
}
\end{center}
The species $C_2^{E}$ in
${\mathbb{Z}_{2}\mbox{-}gpd}^{\mathbb{B}}_1$ sending $x \in
\mathbb{B}$ with $|x| \geq 2$ into
$$C_2^{E}(x)= -2\overline{\mathbb{Z}}_{|x| +1} \sqcup \bigsqcup_{\pi \in Par_d^2(x), d \geq 1}(-2)^{d+1}
\overline{\mathbb{Z}}_{|\pi_d| +1}  \prod_{i=1}^d \prod_{b
\in \pi_i}\overline{\mathbb{Z}}_{|b| +1},$$ satisfies $|C_2^{E}|
=\sum_{n=0}^{\infty}C_{2,n}\frac{x^{n}}{n!}.$ So we obtain a
combinatorial interpretation for the compositional Bernoulli
numbers $C_{2,n}$:
$$C_{N,n}=|C_2^{E}([n])|=\frac{-2}{n+1} + \sum_{\pi \in Par_d^2([n]), d \geq 1}
\frac{(-2)^{d+1}}{(|\pi_d|+1)\prod_{i=1}^d \prod_{b
\in \pi_i}(|b|+1)}.$$

\begin{cor}{\em  $C_{N,1}=1$ and for $n \geq 2$ we have that
$$C_{N,n}= -\frac{-2}{n+1} + \sum_{a_1,...,a_d, d \geq 1}
\frac{(-1)^{d+1}}{(l(a_d)+1)!} \prod_{i=1}^d \prod_{j=1}^{l(a_i)}
\frac{l(a_{i-1})!}{(a_{i1}+1)!....(a_{il(a_i)}+1)!},$$
where $a_1 \vdash l(a_{0}):=n,$ $a_i \vdash l(a_{i-1})$ for $2
\leq i \leq d,$ $l(a_i) \geq 2$ and $a_{ij} \geq 2$.}
\end{cor}

\section{Compositional Bernoulli polynomials}\label{ccbp}

Due to the non-commutativity of composition, there are two natural
compositional generalizations for Bernoulli polynomials, namely,
$$\left(\sum_{n=1}^{\infty}C_{N,n}^{f}\frac{y^{n}}{n!}\right)\circ
f(xy) \mbox{  \ \   and \ \   } f(xy)\circ\left(x,
\sum_{n=1}^{\infty}C_{N,n}^{f}\frac{y^{n}}{n!}\right).$$
The first generalization is easily studied with the help of the
identities:
\begin{eqnarray*}
\left(\sum_{n=1}^{\infty}C_{N,n}^{f}\frac{y^{n}}{n!}\right)\circ f(xy)
&=&\left(\sum_{n=1}^{\infty}C_{N,n}^{f}\frac{y^{n}}{n!}\right)\circ
\left(\sum_{n=1}^{\infty}(f_{n}x^{n})\frac{y^{n}}{n!}\right)\\
&=&\sum_{n=1}^{\infty}\left(\sum_{a \vdash n}
\frac{1}{k!}\left(
\begin{array}{c}
n \\
a_{1},\dots,a_{k}
\end{array}
\right)C_{N,k}^{f}\ f_{a_{1}}x^{a_{1}}\dots f_{a_{k}}x^{a_{k}}
\right)\frac{y^{n}}{n!}\\
&=&\sum_{n=1}^{\infty}\left(\sum_{ a \vdash n}
\frac{1}{k!}\left(
\begin{array}{c}
n \\
a_{1},\dots,a_{k}
\end{array}
\right)C_{N,k}^{f}f_{a_{1}}f_{a_{2}}\dots
f_{a_{k}}x^{n}\right)\frac{y^{n}}{n!}.
\end{eqnarray*}
The second generalization  is less straightforward and is
formalized in our next definition.

\begin{defn}{\em
Let $f=\sum_{n=1}^{\infty}f_{n}\frac{x^n}{n!}\in\mathbb{Q}[[x]]$
and $N \geq 1$ be such that $f_{N}=1$. Compositional Bernoulli
polynomials $C_{N,n}^{f}(x)$ are such that
\begin{equation*}
{\sum_{n=1}^{\infty}C_{N,n}^{f}(x)\frac{y^{n}}{n!}=
f(xy)\circ\left(x,
\sum_{n=1}^{\infty}C_{N,n}^{f}\frac{y^{n}}{n!}\right)}.
\end{equation*}}
\end{defn}

\begin{thm}
{\em Let $f\in\mathbb{Q}[[x]]$ and $N \geq 1$ be such that
$f_{N}=1,$ then }
\begin{equation*}
C_{N,n}^{f}(x)=\sum_{a \vdash n}
\frac{1}{k!}\left(
\begin{array}{c}
n \\
a_{1},\dots,a_{k}
\end{array}
\right)f_{k}C_{N,a_{1}}^{f}\dots C_{N,a_{k}}^{f}x^{k}.
\end{equation*}
\end{thm}
\begin{proof}
\begin{eqnarray*}
\sum_{n=1}^{\infty}C_{N,n}^{f}(x)\frac{y^{n}}{n!}
&=&\left(\sum_{n=1}^{\infty}(f_{n}x^{n})\frac{y^{n}}{n!}\right)\circ
\left(x, \sum_{n=1}^{\infty}C_{N,n}^{f}\frac{y^{n}}{n!}\right)\\
&=&\sum_{n=1}^{\infty}\sum_{a \vdash n}
\frac{1}{k!}\left(
\begin{array}{c}
n \\
a_{1},\dots,a_{k}
\end{array}
\right)f_{k}C_{N,a_{1}}^{f}\dots
C_{N,a_{k}}^{f}x^{k}\frac{y^{n}}{n!}.
\end{eqnarray*}
\end{proof}
We  display compositional Bernoulli polynomials
$C_{N,n}(x)=C_{N,n}^{f}(x)$ for $f(x)=e^{x}$ and $N=1$:
\begin{eqnarray*}
C_{1,0}(x)&=&1, C_{1,1}(x)=x, C_{1,2}(x)=\frac{x^{2}}{2}-x, C_{1,3}(x)=\frac{x^{3}}{6}-x^{2}-2x,\\
C_{1,4}(x)&=&\frac{x^{4}}{24}-\frac{1}{2}x^{3}+\frac{5}{2}x^{2}-6x,\\
C_{1,5}(x)&=&\frac{1}{120}x^{5}-\frac{1}{6}x^{4}+\frac{3}{2}x^{3}-8x^{2}+24x,\\
C_{1,6}(x)&=&\frac{1}{720}x^{6}-\frac{1}{24}x^{5}-\frac{7}{12}x^{4}-
\frac{31}{6}x^{3}+32x^{2}-120x,\\
C_{1,7}(x)&=&\frac{x^{7}}{5040}-\frac{1}{120}x^{6}+\frac{1}{6}x^{5}-\frac{13}{6}x^{4}
+21x^{3}-156x^{2}+720x.
\end{eqnarray*}

\subsection* {Acknowledgment}

Thanks to Jaime Camacaro, Federico Hernandez, Takashi Kimura, Eddy
Pariguan and Domingo Quiroz.

\noindent  hectorblandin@usb.ve\\
Departamento de Matem\'aticas Puras y Aplicadas, Universidad
Sim\'on Bol\'ivar, Caracas 1080-A, Venezuela.\\

\noindent ragadiaz@gmail.com\\
\noindent Grupo de F\'isica-Matem\'atica, Universidad Experimental Polit\'ecnica de las Fuerzas Armadas,
Caracas 1010, Venezuela. \\

\end{document}